\documentclass[a4paper,11pt]{article}
\usepackage{graphicx}
\usepackage{geometry}
\usepackage{booktabs}
\usepackage{titlesec}
\usepackage{setspace}
\usepackage{natbib}
\usepackage{helvet} 
\usepackage{float}
\usepackage[T1]{fontenc}
\usepackage[utf8]{inputenc}
\usepackage[french]{babel}
\usepackage{csquotes}

\usepackage{amsmath}   
\usepackage{amssymb}   

\author{}
\date{}

\begin{document}


{\LARGE{\begin{center}{\textbf{Uncovering commuting flows in Bike Sharing Systems}}\end{center}}}

\vspace{1em}

{\centering \textbf{Mohamadou Salifou$^{1,2*}$}\par}

{\centering $^{1}$Univ. de Toulouse, Univ. Toulouse 2, CNRS, LISST, 5 Allées Antonio Machado, TOULOUSE CEDEX 9, Toulouse 31058, France. \par}
{\centering $^{2}$Univ. de Toulouse, F\'ed\'eration ENAC ISAE SUPAERO ONERA, 7 avenue Edouard Belin, Toulouse Cedex 4, Toulouse 31055, France. \par}
{\centering *mohamadou.salifou@univ-tlse2.fr \par}

\section*{Introduction}

Since the 2000s, Bike Sharing Systems (BSS) have experienced rapid expansion in cities around the world. They are now a central component of sustainable mobility policies, providing a flexible and active travel solution that complements public transportation. Among the different BSS models, one way systems ($O \to D$) in which users can borrow a bicycle from one station and return it to another are the most widespread \citep{gu2019or,de2019contradictions}.\\

This configuration generates complex flows between stations, whose analysis has been made possible by the availability of large scale datasets describing each trip individually. A significant portion of the literature has therefore focused on the global determinants of BSS usage, examining the spatio-temporal variations in demand, as well as the influence of weather conditions, the urban environment, and the socio-demographic characteristics of neighborhoods \citep{gu2019or,de2019contradictions,todd2021global}. However, these aggregated approaches often overlook the fine structure of flows between station pairs, and in particular the dynamics of loop trips, even though these represent a valuable indicator of individual behaviors and their regularities \citep{nello2021scaling,wilkesmann2023determinants}. Although they represent a relatively small share of the total volume (ranging from $5\%$ to $20\%$ depending on the city), such trips are far from anecdotal \citep{bordagaray2016capturing,zhang2018mining,kou2020quantifying}. They are often associated with recreational or tourist uses, but several studies have shown that they may also correspond to commuting or repeated utilitarian trips, particularly among regular subscribers \citep{o2014mining,bordagaray2016capturing,zhang2018mining,lee2025measuring}.\\

This distinction is, however, difficult to establish in $O \to D$ systems, because the absence of individual identifiers prevents confirming with certainty that an outward trip and a return trip are performed by the same user. This methodological difficulty explains why some studies simply exclude loop trips \citep{ma2019spatiotemporal,kou2020quantifying}, while others attempt to infer their meaning using spatial or contextual proxies such as the environmental characteristics of stations \citep{zhao2015exploring,liu2019associations}. As a result, the understanding of round-trip behaviors in BSS remains incomplete.\\

In recent work, \citep{detourjouve2025} showed that modeling trip durations using mixture models, particularly log-normal mixtures, makes it possible to highlight behavioral heterogeneity within station pairs. By explicitly distinguishing the two directions $O \to D$ and $D \to O$, one observes not only homogeneous or heterogeneous station pairs, but also pairs exhibiting directional asymmetry in the very structure of trip duration distributions.\\

From the same perspective, how can the heterogeneity of trip duration distributions observed jointly in both directions between a pair of stations be exploited to identify temporally organized mobility structures ? More specifically, can the dominant components of mixture models be used to reveal regular and constrained usage patterns, compatible with commuting behaviors, without relying on individual or survey data?

To address this question, we propose an approach based on the joint exploitation of temporal and directional dimensions to identify station pairs likely to support commuting flows. This approach makes it possible to connect the statistical analysis of duration distributions with a behavioral interpretation of observed flows, without requiring individual-level data.\\

In a distinct framework, the second part of the article introduces a probabilistic model explicitly linking morning and evening flows. By exploiting the structural dependence between these two periods using a Hurdle Poisson - Binomial framework \citep{bayart2012interet, feng2021comparison} , we propose a direct estimation method for the proportion of morning trips that result in a return trip on the same day. This approach, based on aggregated flows rather than mixture model components, provides a complementary tool for quantifying the relative importance of commuting movements and for exploring predictive applications, particularly in the operational management of BSS.

Once this proportion has been estimated, it can be used for predictive purposes. More specifically, knowledge of the proportions of morning trips ($O \to D$, $O \in V$), that lead to an effective return at the end of the day makes it possible to predict the number of bicycle uses expected in the late afternoon departing from station $D$. Thus, the information extracted from the flows observed in the morning is used to anticipate future demand across the network, station by station. One of the major advantages of this approach lies in the possibility of direct empirical validation : the predicted volumes of trips in the second part of the day can be compared with the actual observations in the data, thereby allowing the relevance and robustness of the proposed model to be assessed.

\section{Methodology}

\subsection{Data description and preprocessing}
The data used in this study come from the BSS of Toulouse for the year 2022. They include the departure and arrival stations, the start and end durations of each trip, as well as the corresponding date. Of the 3 983 301 trips initially recorded, we excluded those with a duration exceeding thirty minutes (121 522 trips), in order to eliminate potentially atypical uses, such as bikes left stationary for long periods instead of being actively used. We also removed trips where the origin and destination stations were identical (136 809 trips, representing $3.43\%$), due to the absence of individual level data, which prevents us from determining whether these were round trips made by the same user or two distinct trips by different users. After this filtering process, the final dataset includes 3 208 952 valid trips.\\

Each station pair ($O, D$) is considered as a unit of analysis. In addition to the variables available for individual trips, we created the following: direction of trip ($O \to D$ ou $D \to O$), day type (weekday or weekend), and time of day (night, morning, midday, evening). To establish a consistent and observation independent trip direction, stations were ordered according to their distance from the city center, represented by Place du Capitole. For each station pair, the station closest to the Capitole is systematically designated as station $O$, while the farthest station is designated as station $D$. This convention allows a unique geographic orientation for each pair and makes it possible to distinguish the trip directions. Trips are then classified as $O \to D$ or $D \to O$ depending on whether the departure station corresponds to station $O$ or station $D$, respectively.
We retained only station pairs with at least 100 trips in each direction over the year. In total, 9 176 round-trip pairs were included in the analyses.


\subsection{Statistical signatures of commuting trips}

We model the distribution of trip durations in each direction using a log-normal mixture model and retain the station pairs for which heterogeneous behaviors are observed in both directions.  As shown by \cite{detourjouve2025}, the first component of the mixture model typically corresponds to the reference itinerary, which is consistent with the fastest route suggested by OpenStreetMap. For each observed trip, we therefore compute the posterior probability of belonging to the mixture component associated with this typical route.\\

These posterior probabilities are then modeled as a function of contextual variables using a logistic regression framework:

$$
\Pr(Z = 1)
=
\beta_0
+ \beta_1 \text{Morning}
+ \beta_2 \text{Midday}
+ \beta_3 \text{Evening}
+ \beta_4 \text{Weekend}
+ \beta_5 \text{Direction}
+ \text{interaction terms}.
$$

where $Z = 1$ denotes membership in the short (or fast) trip component.\\

A pair of stations is considered a potential candidate for commuting flows if it simultaneously satisfies the following statistical criteria.

\paragraph{Strong morning and evening effects on weekdays :}

The coefficients associated with peak commuting hours satisfy

\[
\beta_{\text{Morning}} > 0
\quad \text{and} \quad
\beta_{\text{Evening}} > 0,
\]

and are statistically significant. In contrast, the nighttime period (taken as the reference category) is associated with a low probability of belonging to the fast trip component.

\paragraph{Negative or negligible weekend effect :}

The coefficient associated with weekends satisfies

\[
\beta_{\text{Weekend}} < 0,
\]

and is statistically significant. The reduction of such trips during weekends is consistent with utilitarian mobility patterns rather than recreational usage.

\paragraph{Directional asymmetry between morning and evening :}

The most discriminative feature of commuting flows lies in the interaction between time of day and trip direction. In particular, we consider interaction terms of the form

$$
\text{Morning} \times \text{Direction},
\qquad
\text{Evening} \times \text{Direction}.
$$
For typical commuting patterns, we expect

\[
\beta_{\text{Morning}\times\text{Direction}} > 0,
\qquad
\beta_{\text{Evening}\times\text{Direction}} < 0.
\]

This pattern reflects the characteristic temporal reversibility of commuting flows, with flows concentrated in one direction during the morning and in the opposite direction during the evening.

\paragraph{Weakly structured midday period :}

Finally, the coefficient associated with the midday period typically satisfies

\[
\beta_{\text{Midday}} \approx 0
\]

or is not statistically significant. This indicates that midday trips tend to be more heterogeneous and less structured than peak commuting trips.

\subsection{Morning and evening flow model}

We propose a model that makes it possible to estimate the proportion of morning trips that lead to a return trip in the evening in a BSS. In addition, we sometimes observe no activity during the morning, resulting in so called structural zeros. These may correspond to public holidays, school vacations, exceptional weather conditions, or temporary station closures (inactive stations). Such zeros don't reflect mobility behavior but rather the absence of the underlying mobility process. To model this phenomenon, we use a \textit{Poisson Hurdle} model \citep{bayart2012interet, feng2021comparison}.

\subsubsection{Morning flow model}

We introduce an indicator variable

$$
Z_i =
\begin{cases}
1 & \text{if the station pair is active on day } i, \\
0 & \text{otherwise (structural zero).}
\end{cases}
$$

We assume

$$
P(Z_i = 1) = \pi.
$$

Thus, $\pi$ represents the probability that a day is active for the station pair.\\

The variable $N_i$, representing the number of morning trips between stations $O$ and $D$, is defined as

$$
N_i =
\begin{cases}
0 & \text{if } Z_i = 0,  \text{with probability } 1-\pi\\
Y_i & \text{if } Z_i = 1,  \text{with probability } \pi
\end{cases}
$$

where

\[
Y_i \sim \text{Poisson}(\lambda) \text{ truncated at } 0.
\]

Thus, $N_i$ follows a Poisson Hurdle model.\\

For a zero-truncated Poisson distribution, we have

\begin{equation*}
\begin{aligned}
\mathbb{E}[Y_i] &= \frac{\lambda}{1-e^{-\lambda}}, \\
Var(Y_i) &= \frac{\lambda}{1-e^{-\lambda}}\left(1 - \frac{\lambda e^{-\lambda}}{1-e^{-\lambda}}\right).
\end{aligned}
\end{equation*}

Consequently,

$$
\mathbb{E}[N_i] = \pi \frac{\lambda}{1-e^{-\lambda}}.
$$

\begin{equation*}
Var(N_i) = \pi \frac{\lambda}{1-e^{-\lambda}} \left(1 - \frac{\lambda e^{-\lambda}}{1-e^{-\lambda}}\right)
+ \pi(1-\pi) \left(\frac{\lambda}{1-e^{-\lambda}}\right)^2.
\end{equation*}

The estimates of $\pi$ and $\lambda$ are

$$
\hat \pi = \frac{\text{number of days with } N_i>0}{\text{total number of days}} = \frac{n_+}{n}.
$$

For active days, we compute the empirical mean (parameter of the truncated Poisson)

$$
\bar N_+ = \frac{1}{n_+} \sum_{N_i>0} N_i
$$

The parameter $\lambda$ is the solution of

$$
\bar N_+ = \frac{\lambda}{1-e^{-\lambda}},
$$

which is solved numerically.

\subsubsection{Evening flow model}

Among morning trips, a proportion $p$ corresponds to individuals returning in the evening:

$$
R_i \mid N_i \sim \text{Binomial}(N_i,p).
$$

Trips independent of the morning flow may also occur in the evening:

$$
X_i \sim \text{Poisson}(\mu), \quad X_i \perp N_i.
$$

The total evening flow is therefore

$$
M_i = R_i + X_i.
$$

\subsubsection{Estimator of $p$}

By linearity of covariance,

$$
Cov(N_i,M_i) = Cov(N_i,R_i) + Cov(N_i,X_i).
$$

Since $X_i$ is independent of $N_i$, we obtain

$$
Cov(N_i,M_i) = Cov(N_i,R_i).
$$

Conditionally on $N_i$, we have

$$
R_i \mid N_i \sim \text{Binomial}(N_i,p),
$$

which implies

$$
\mathbb{E}[R_i \mid N_i] = pN_i.
$$

Using the law of total covariance,

$$
Cov(N_i,R_i)
=
Cov\bigl(N_i,\mathbb{E}[R_i\mid N_i]\bigr)
+
\mathbb{E}\bigl[Cov(N_i,R_i\mid N_i)\bigr].
$$

Since $N_i$ is constant conditionally on itself,

$$
Cov(N_i,R_i\mid N_i)=0.
$$

Therefore,

$$
Cov(N_i,R_i)
=
Cov(N_i,pN_i)
=
p\,Var(N_i).
$$

Consequently,

$$
Cov(N_i,M_i) = p\,Var(N_i).
$$

Hence, a method of moments estimator for the parameter $p$ is

$$
\hat p =
\frac{\widehat{Cov}(N_i,M_i)}
{\widehat{Var}(N_i)}.
$$

More explicitly,

\begin{equation*}
\hat p = \frac{\widehat{Cov}(N_i,M_i)}
{\pi \frac{\lambda}{1-e^{-\lambda}} \left(1 - \frac{\lambda e^{-\lambda}}{1-e^{-\lambda}}\right)
+ \pi(1-\pi) \left(\frac{\lambda}{1-e^{-\lambda}}\right)^2}.
\end{equation*}

Once $\hat p$ is estimated, the evening flow independent of the morning is given by

$$
\hat \mu = \bar M - \hat p \, \bar N.
$$

The estimator $\hat p$ therefore remains interpretable as the proportion of morning trips that lead to an evening return trip, even in the presence of days with no observed flow.

\section{Results}

\subsection{Overview of the regression analysis}

\subsubsection{Global statistical patterns}

Across the entire system, several systematic temporal patterns emerge.

Table~\ref{tab:global_patterns} reports the proportion of station pairs exhibiting significant effects for each explanatory variable.

\begin{table}[h]
\centering
\caption{Proportion of station pairs exhibiting significant temporal effects.}
\begin{tabular}{lcc}
\hline
Variable & Number of pairs & Percentage \\
\hline
Morning peak effect & 4784 & 52.14\% \\
Evening peak effect & 4089 & 44.56\% \\
Negative weekend effect & 5984 & 65.21\% \\
Directional asymmetry & 3457 & 37.67\% \\
\hline
\end{tabular}
\label{tab:global_patterns}
\end{table}

These results indicate that peak-hour effects are widespread across the system, while strong directional asymmetry remains more localized.

\subsubsection{Typology of station pairs}

Based on the regression results, station pairs can be grouped into four main categories.

\begin{table}[h]
\centering
\caption{Typology of station pairs based on regression signatures.}
\begin{tabular}{lcc}
\hline
Type of pair & Criteria & Number of pairs \\
\hline
Strong commuting flows & all criteria satisfied & 1180 \\
Probable commuting & morning + evening peaks & 2050 \\
Leisure oriented trips & positive weekend effect & 1700 \\
Diffuse mobility & no clear pattern & 4246 \\
\hline
\end{tabular}
\end{table}

This classification highlights the coexistence of structured commuting flows and more flexible leisure mobility within the system. We observe that 1 185 out of 9 176 pairs (approximately 12.91$\%$) exhibit a strong pattern of commuter flows. When the criteria are relaxed, this number increases to 3 238 pairs (about 35.28$\%$). Trips that are more related to leisure activities characterized by a predominance of weekend trips and an absence of structured morning and evening usage account for 18.54$\%$ of the pairs.

\subsubsection{Spatial structure of commuting flows}

\begin{figure}[H]
    \centering
    \includegraphics[width=\textwidth]{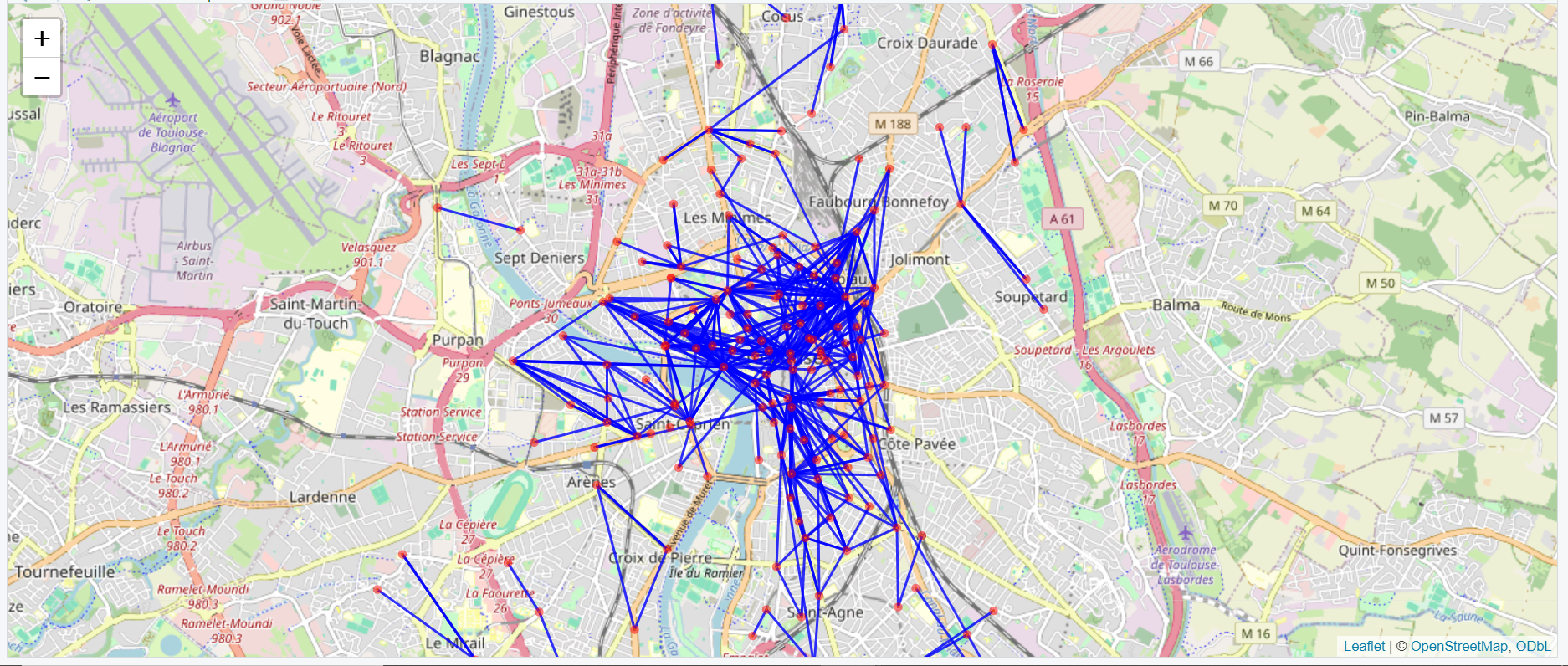}
    \caption{Spatial structure of commuting flows. The red points correspond to the stations of the BSS, while the blue segments connect pairs of stations between which trips have been observed.}
    \label{fig:temporal_profile}
\end{figure}

Mapping the identified commuting pairs reveals a clear spatial organization (see figure \ref{fig:temporal_profile}).

The figure represents the spatial structure of station pairs identified as supporting commuting flows. A strong concentration of connections can be observed in the central area of the city, indicating that this zone acts as a major mobility hub. This suggests that many users living in the city center use the system to travel to their places of activity. Indeed, several stations located in the central districts appear to be highly connected, which indicates that they function as hubs in the network, concentrating a large number of departures and arrivals.

Moreover, many segments link peripheral neighborhoods to the city center, forming a globally radial structure that is typical of daily urban mobility patterns. These spatial corridors suggest the presence of regular flows between residential areas located on the outskirts of the city and activity or employment zones situated in the central districts. It can also be observed that some connections extend along major urban axes, which may reflect the influence of the road network configuration or the location of major activity centers. Overall, the figure highlights a strongly centralized spatial organization of the BBS, consistent with commuting mobility patterns in which users travel from residential areas to central districts during working periods and return in the opposite direction at the end of the day.

\subsection{Estimated proportion of round trips}

\begin{figure}[H]
    \centering
    \includegraphics[width=\textwidth]{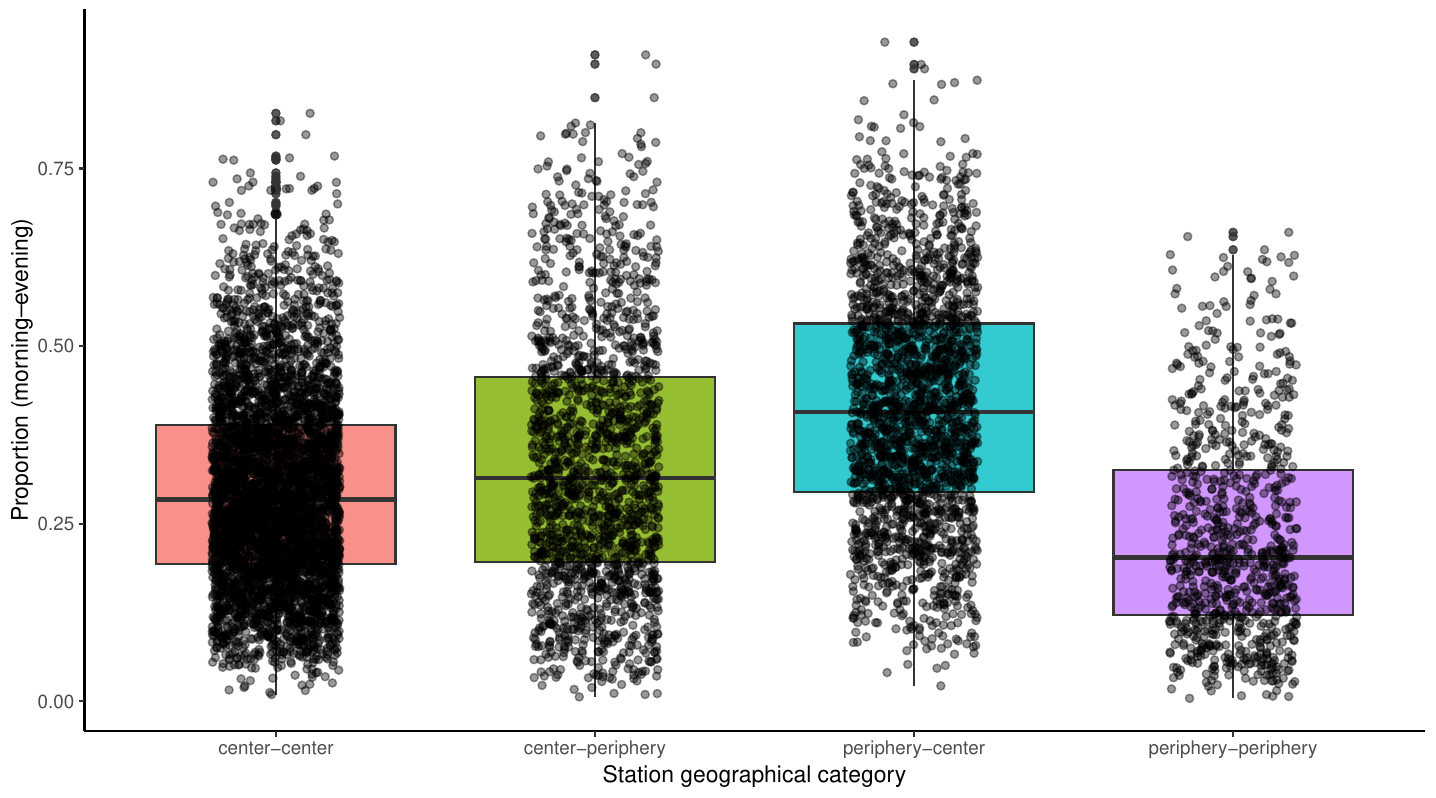}
    \caption{Distribution of the proportion $p$ of morning trips that result in a return trip in the evening, according to the geographical category of the origin and destination stations.}
    \label{fig:geography}
\end{figure}

Clear differences appear across the different types of trips. Trips from the periphery to the city center exhibit the highest median values of $p$, indicating that a large share of morning trips in this configuration is followed by a return trip in the evening. This pattern is consistent with a commuting dynamic in which cyclists travel toward central urban areas in the morning and return toward peripheral areas at the end of the day. Trips from the center to the periphery and trips within the center show intermediate proportions, suggesting more heterogeneous mobility behaviors where morning trips don't systematically correspond to a daily round trip. In contrast, trips between peripheral stations display the lowest values of $p$, indicating that morning trips in these areas are less frequently associated with an evening return. Overall, these results highlight a spatial structure in cycling mobility in which trips involving the city center are more strongly associated with daily round trip commuting patterns.

\section{Conclusion and recommendations}

In conclusion, this study proposes two methodological approaches. First, it uses statistical patterns in trip durations distributions to identify station pairs for which the fastest trip component likely reflects commuting behavior. Second, it introduces an estimator of the proportion of morning trips that result in a return trip in the evening, even in the absence of individual level data that would allow users to be directly identified.\\ 

The results provide several implications for the operational management of BBS. The analysis of the posterior probability of belonging to the component associated with the fastest route durations shows that these trips occur more frequently at specific times of the day and depending on the type of day. This suggests that the fastest component can reasonably be interpreted as representing utilitarian trips. These findings can help operators anticipate spatial and temporal imbalances in bicycle availability across stations. For instance, stations located in residential areas may experience high demand during morning peak hours, while stations near employment centers may become more congested upon users’ arrival. Anticipating these directional flows could therefore improve bicycle redistribution strategies.

The proposed probabilistic framework allows the inference of mobility patterns even when individual trip identifiers are unavailable. This approach could be applied to other bike-sharing systems with similar aggregated data.

\section*{Acknowledgements}

This research was supported by Région Occitanie (n° DOMSUB23003131) and Université de Toulouse (n°155 2024).

\bibliographystyle{plain}
\bibliography{refs}

\end{document}